\documentclass[submission]{pmsart}[2006/01/13]

\allowdisplaybreaks[4]
\usepackage{color}
\usepackage[colorlinks=true, citecolor=blue, linkcolor=blue]{hyperref}

\usepackage{xparse}
\NewDocumentCommand{\ceil}{s O{} m}{%
  \IfBooleanTF{#1} 
    {\left\lceil#3\right\rceil} 
    {#2\lceil#3#2\rceil} 
}
\NewDocumentCommand{\floor}{s O{} m}{%
  \IfBooleanTF{#1} 
    {\left\lfloor#3\right\rfloor}
    {#2\lfloor#3#2\rfloor}
}

\definecolor{c20}{rgb}{0.,0.7,0.}
\definecolor{c30}{rgb}{0.,0.,1.}
\definecolor{c40}{rgb}{1,0.1,0.7}
\definecolor{c50}{rgb}{1,0,0}
\definecolor{c60}{rgb}{1,0.9,0.1}

\def\peH#1{\textcolor{c20}{#1}}

\def\ii#1{\textcolor{c50}{#1}}

\def\EHC#1{\textcolor{black}{#1}}
\def\KD#1{\textcolor{c30}{#1}}

\def\KD#1{#1}

\def\Ed#1{\textcolor{c30}{#1}}
\def\Ee#1{\textcolor{c50}{#1}}

\def\Ee#1{#1}

\def\K#1{\textcolor{cyan}{#1}}
\def\K#1{\textcolor{black}{#1}}
\def\ii#1{\textcolor{black}{#1}}
\def\KK#1{\textcolor{cyan}{#1}}
\def\KK#1{\textcolor{black}{#1}}

\def\Ed#1{#1}

\def\peH#1{#1}

\newcommand{\ve}{\varepsilon}

\newcommand{\abs}[1]{\left\lvert #1 \right\rvert}

\newcommand{\E}[1]{\mathbb{E}\left\{#1\right\}}
\newcommand{\EE}[1]{\mathbb{E}\{#1\}}

\newcommand{\pk}[1]{\mathbb{P} \left\{ #1 \right\} }

\newcommand{\R}{\mathbb{R}}
\newcommand{\Z}{\mathbb{Z}}

\newcommand{\N}{\mathbb{N}}
\newcommand{\inr}{\in \R}
\newcommand{\inn}{\in \N}
\newcommand{\ldot}{,\ldots,}

\newcommand{\limit}[1]{\lim_{#1 \to   \infty}}

\newcommand{\BQN}{\begin{eqnarray}}
\newcommand{\EQN}{\end{eqnarray}}
\newcommand{\BQNY}{\begin{eqnarray*}}
\newcommand{\EQNY}{\end{eqnarray*}}
\def\bqny#1{ { \begin{eqnarray*} #1 \end{eqnarray*}}}
\def\bqn#1{ { \begin{eqnarray} #1 \end{eqnarray}}}

\newcommand{\BS}{\begin{sat}}
\newcommand{\ES}{\end{sat}}
\newcommand{\BT}{\begin{theo}}
\newcommand{\ET}{\end{theo}}
\newcommand{\BK}{\begin{korr}}
\newcommand{\EK}{\end{korr}}
\newcommand{\EQD}{\stackrel{d}{=}}

\newcommand{\BD}{\begin{de}}
\newcommand{\ED}{\end{de}}

\newcommand{\BIT}{\begin{itemize}}
\newcommand{\EIT}{\end{itemize}}
\newcommand{\BDI}{\begin{description}}
\newcommand{\EDI}{\end{description}}

\newcommand{\BRM}{\begin{remarks}}
\newcommand{\ERM}{\end{remarks}}

\newcommand{\BEL}{\begin{lem}}
\newcommand{\EEL}{\end{lem}}

\newcommand{\nelem}[1]{{Lemma \ref{#1}}}
\newcommand{\neprop}[1]{{Proposition \ref{#1}}}
\newcommand{\netheo}[1]{{Theorem \ref{#1}}}

\newcommand{\prooftheo}[1]{ \textsc{\bf Proof of Theorem} \ref{#1}:}
\newcommand{\proofprop}[1]{\textsc{\bf Proof of Proposition} \ref{#1}:}
\newcommand{\prooflem}[1]{\textsc{\bf Proof of Lemma} \ref{#1}:}

\newcommand{\COM}[1]{}

\newcommand{\QED}{\hfill $\Box$}

\def\xiW{\xi_W}

\newtheorem{theo}{theorem}[section]
\newtheorem{Cond}[theo]{Condition}

\def\IF{\infty}

\date{}

\def\Var{\text{Var}}

\def\FRE{\mbox{Fr\'{e}chet }}
 \def\YD#1{Y^\delta{(#1)}}
 \def\xiWD#1{\xi_W^\delta ( #1 ) }
 \def\XD#1{X^\delta ( #1 ) }
 \def\WD#1{W^\delta ( #1 ) }
 \def\XWD#1{X^\delta ( #1 ) }
\def\HWD{\mathcal{H}_W^\delta}
\def\HWDA{\mathcal{H}_W^0}
\def\HWDAA{\mathcal{H}_W}

\def\EID{\theta_W^\delta}
\def\EIDC{\widetilde{\theta_W^\delta}}
\def\EIDD{\widehat{\theta_W^\delta}}

\volume{0}%
\fasc{0}%
\years{0000}%
\pages{000--000}%
\firstpage{1} %
\lastrevision{xx.xx.xxxx}

\title{On Extremal Index  of max-stable stationary processes}
\titlethanks{}
\dedicated{To Tomasz Rolski,\\
thankful for all the support and ideas you shared with us!}
\ShortTitle{On Extremal Index  of Max-stable Stationary Processes} 
\ShortAuthors{K. D\c{e}bicki and E. Hashorva}
\NumberOfAuthors{2}

\received{24.11.2016}
\FirstAuthor{Krzysztof D\c{e}bicki}
\FirstAuthorCity{Wroc\l aw}
\FirstAuthorAffiliation{Mathematical Institute, University of Wroc\l aw}
\FirstAuthorAddress{pl. Grunwaldzki 2/4,\\ 50-384 Wroc\l aw, Poland}
\FirstAuthorEmail{Krzysztof.Debicki@math.uni.wroc.pl}                              
\FirstAuthorThanks{}                             

\SecondAuthor{Enkelejd  Hashorva} 
\SecondAuthorCity{Lausanne} %
\SecondAuthorAffiliation{University of Lausanne}%
\SecondAuthorAddress{B\^{a}timent Extranef, UNIL-Dorigny,\\ 1015 Lausanne, Switzerland}%
\SecondAuthorEmail{Enkelejd.Hashorva@unil.ch}%
\SecondAuthorThanks{}%

\ThirdAuthor{}
\ThirdAuthorCity{}%
\ThirdAuthorAffiliation{}%
\ThirdAuthorAddress{}%
\ThirdAuthorEmail{}
\ThirdAuthorThanks{}%

\FourthAuthor{}
\FourthAuthorCity{}%
\FourthAuthorAffiliation{}%
\FourthAuthorAddress{}%
\FourthAuthorEmail{}%
\FourthAuthorThanks{}%

\FifthAuthor{} 
\FifthAuthorCity{} %
\FifthAuthorAffiliation{}%
\FifthAuthorAddress{} %
\FifthAuthorEmail{}%
\FifthAuthorThanks{}%


\MSCcodes{Primary: 60G15; Secondary: 60G70.}               

\begin{document}
\maketitle
\begin{abstract}
In this contribution we discuss the relation between Pickands-type constants defined for \Ee{certain}
Brown-Resnick stationary process  $W(t),t\inr$ as
\[\HWD=
\limit{T} T^{-1}  \E{ \sup_{t\in \delta \Z \cap  [0,T]} e^{W(t)} }, \ \delta \ge 0\]
\ii{(set $0 \Z=\R$ if $\delta=0$)} and the extremal index of the associated max-stable stationary process
 $\xiW$.
 We derive several new formulas and obtain lower bounds for $\HWD$ if $W$
 is a Gaussian or a L\'evy process. As a by-product we show an
 interesting relation between Pickands constants
 and lower tail probabilities for fractional Brownian motions.
\end{abstract}

\section{Introduction}
The motivation for this contribution comes from the importance and the intriguing properties of the classical Pickands
constants $\HWD$, which are defined for any $\delta \ge 0$ by (interpret  $0 \Z$ as $\R)$
\BQN \label{defP}
\HWD = \limit{T} \frac{1}{T} 
\E{ \sup_{t\in \delta \Z \cap  [0,T]} e^{W(t)} } ,
  \EQN
where  $$W(t)= \sqrt{2}B_\alpha(t)- \abs{t}^\alpha, \quad t\inr,$$
 with $B_\alpha$ a standard fractional
Brownian motion with Hurst index $\alpha \in (0,2]$,
\K{that is a centered Gaussian process with stationary increments and variance function $\Var(B_\alpha(t))=|t|^\alpha, t\inr$.}

It is well-known (but not trivial to prove) that
$\HWD$  is finite  and positive for any $\delta\ge 0$. The only values known for $\HWD$ are for  $\delta=0$ and $\alpha \in \{1, 2\}$, see e.g., \cite{PickandsB,Pit20}.  Suprisingly, Pickands and related constant appear \ii{in numerous unrelated asymptotic problems,} see e.g., the recent papers \cite{Delrome, Harper1, Harper2, KrzKos}. \\
The contribution \cite{DiekerY} derived a new formula for Pickands constants, which in fact indicates
a direct connection between those contants and max-stable stationary processes, see \cite{SBK}. The definition of $\HWD$ in \eqref{defP} is extended in \cite{SBK}
for some general process $W$, provided that it defines a max-stable and stationary process.  More precisely, assume throughout in the sequel that
\BQN \label{Ws}
 W(t)= B(t) - \ln   \E{e^{B(t)}},  \quad t\inr,
\EQN
where $B(t), t\inr$ is a random process on the space $D$ of c\`adl\`ag functions \\
$f:\R\to\R$ with
\BQN \label{B0}
B(0) = 0, \quad  \E{e^{B(t)}}< \IF,  \quad  t\inr.
\EQN
Hence $X(t)=e^{W(t)}$ satisfies $X(0)=1$ almost surely, and $\E{X(t)}=1, t\inr$.\\
 If $\Pi=\sum_{i=1}^\IF \ve_{P_i}$ is a Poisson point process (PPP) with intensity $x^{-2} dx$ on $(0,\IF)$, and $X_i=e^{W_i},i\ge 1$ are independent  copies of the random \Ee{process} $X=e^W$ being independent of $\Pi$, then the random process $\xiW$ defined by
\BQN\label{eA}
 \xiW (t)= \max_{i\ge 1} P_i X_i(t)=\max_{i\ge 1} P_i e^{ W_i(t)}  , \quad t\inr
\EQN
has unit \FRE marginals and  is max-stable.
Here 
$\ve_x$ denotes the unit Dirac measure at $x\in\R$.\\
Adopting the definition in \cite{kab2009}, we shall refer to $W$ as the
{\it Brown-Resnick stationary process}
whenever the associated max-stable process $\xiW$ is stationary.
Note that stationarity of $\xiW$ means that $\{ \xiW (t),t\inr\}$ and $\{ \xiW (t+h),t\inr\}$ have the same distribution for any $h\inr$. 

\COM{holds for any $\delta\ge 0$, where $\HWD \in (0,\IF)$ is the classical Pickands constant .

Let

$ X(t)= e^{W(t)}, t\in \R$ be \Ee{a} random process  such that $\E{X(t)}=1, t\inr$.
In the special case that $W(t)= \sqrt{2}B_\alpha(t)- \abs{t}^\alpha$ with $B_\alpha$ a standard fractional
Brownian motion with Hurst index $\alpha \in (0,2]$,
it is well-known that
\BQN \label{defP}
H_{W}^\delta([0,T]):=  \frac{1}{T} \E{ \sup_{t\in \delta \Z \cap  [0,T]} e^{W(t)} } \to \HWD, \quad T \to \IF.
  \EQN
holds for any $\delta\ge 0$, where $\HWD \in (0,\IF)$ is the classical Pickands constant (interpret  $0 \Z$ as $\R)$.
As shown in the recent contribution \cite{SBK} the study of Pickands-type constants $\HWD$, which are
the limit of  the functional $H_{W}^\delta([0,T])$ as $T\to\infty$
is of interest for both the extreme value theory of max-stable processes and
 extremes of Gaussian processes. \\
 More specifically, let in the sequel $W$ be a drifted process defined as follows (we keep the same notation as in  \cite{SBK})
}
In the sequel, for the case $\delta=0$  we shall assume that
$$\E{ \sup_{t\in K} e^{W(t)} }   < \IF$$
  for any compact
$K \subset \R$.    A direct consequence of stationarity of $\xiW$ and
{the fact that
for any $t_1 \ldot t_n \inr$ and $x_1 \ldot x_n >0$, see e.g., \cite{DM,MolchanovSPA}
\BQN\label{base}
\pk{ \xiW(t_i) \le x_i, \forall i\le n} = e^{- \E{\max_{1 \le i\le n} \bigl(e^{W(t_i)}/x_i \bigr)}}
\EQN
}
is that, for any $b\ge 0,\delta \ge 0,  T>0$ we have
\BQNY
 H_{W}^\delta([0,T]) &:=&\E{ \sup_{t\in \delta \Z \cap  [0,T]} e^{W(t)} }
=
 \E{ \sup_{t\in \delta \Z \cap  [b,b+T]} e^{W(t)} }.
\EQNY
Consequently, $\mathcal{H}_{W}^\delta$ defined in \eqref{defP} exists and  is given by (see \cite{SBK})
\BQN
 \label{fino}
 \HWD =	 \inf_{T>0} \frac{1}{T} H_{W}^\delta([0,T]) \in \Ee{[}0,\IF).
 \EQN
Note that if $\delta>0$, then \eqref{fino} \Ee{implies that
	$$\HWD \le  \frac {H_{W}^\delta([0,\delta- \ve])}{ \delta- \ve}= \frac 1 {\delta- \ve} $$
	 for any $\ve \in (0,\delta)$, hence letting $\ve$ tend to 0 yields} $\HWD \in [0, \delta]$.

\ii{Interestingly}, $\HWD$ is related to the extremal index of the stationary process
$$\xiWD{t}= \xiW(\delta t), \quad t \in \Z, \delta >0,$$
where we set  $\xiWD{t}=\xiW(t)$ if $\delta=0$. Indeed, \K{by (\ref{base})}
\begin{eqnarray} \label{FreLim}
\limit{T} \pk{ \max_{i\in  \delta \Z \cap [0,T]} \xiW(t) \le T x}
&=&
e^{- \limit{T}
 \E{\max_{i\in  \delta \Z \cap [0,T]}  \bigl( e^{W(i)}/T\bigr)} \frac{1}{x}}\nonumber\\
& =&
 \Bigl( e^{- \frac{1}{x}} \Bigr)^{\HWD}  , \quad x>0.
\EQN
\COM{where
$$ \mathcal{H}^\delta_X =\limit{T} \E{\max_{i\in  \delta \Z \cap [0,T]}  \frac{X(i)}{T}}=
\inf_{T> 0}
 \E{\max_{i\in  \delta \Z \cap [0,T]}  \frac{X(i)}{T}}  \in [0, \IF).
 $$
 }
 Thus the \FRE limit result in \eqref{FreLim}, which is already shown in
\cite{stoev2010max} (see also \cite{MIKYuw}[Proposition 3.1] and \cite{DM})
 states that the extremal index of the stationary process
$\xiWD{t} ,t \in \Z$ is given  for any $\delta>0$ by
\BQN \label{FreLim2}
 \EID= \delta \HWD \in [0,1].
 \EQN
Clearly, the constant $\HWD $ is positive if and only if the extremal index $\EID$ of the stationary process $\xiW^\delta$ is positive.

Numerous papers in the literature have discussed the calculation and estimation of extremal index  of stationary processes, see e.g., the recent articles \cite{MR2248577,MIKYuw,MikoshZ, MR3179972, Soulier,Jakub} and the references therein. \\
\peH{The primary goal of this \Ee{contribution is to study Pickands type constants $\HWD$ by exploring the properties of the extremal index $\EID$. In particular, we are interested in establishing tractable conditions that guarantee the positivity of $\HWD$.}\\
By our assumptions it is clear that $\xiW^\delta$ is stationary and {\it jointly regularly varying}, hence in view of  \cite[Theorem 2.1]{BojanS} (see also \cite{MR3301293}), there exists the so-called {\it tail process}
 $$\YD{i} , \quad i\in \Z$$
 of the stationary process $X$, which was introduced in \cite{BojanS}. 
   It turns out that for any $m\le n, m,n\in \Z$ we have the following stochastic representation
\BQN \label{stochRep}
(\YD{m} \ldot \YD{n}) &\EQD & ( \ii{\mathcal{P}}\XWD{m} \ldot  \ii{\mathcal{P}}\XWD{n}),
\EQN
$ X^\delta(i):= e^{W(\delta i)}, i\in \Z,$
with $\ii{\mathcal{P}} $ a unit Pareto random variable with survival function $1/x, x>1$ being independent of the process $X$.\\
Under the {\it finite mean cluster size condition} (see below Condition \ref{C}) and condition $\ii{\mathcal{A}}(a_n)$, see \cite{BojanS,MR3025708,Danijel}, it follows that $\EID$ is positive, see the seminal contribution \cite{BojanS}.\\
 We shall show the positivity of the extremal index under a weaker condition, namely
 \K{supposing that}
\BQN  \label{asC}
\lim_{\abs{z} \to \IF, z\in \Z } W(z \delta ) =  - \IF
\EQN
holds almost surely for $\delta\in (0,\IF)$. In our derivations the next simple result is crucial:
\begin{lemma} \label{propAA}If $r_n,n\ge 1$ are positive integers such that
\[\limit{n} r_n= \limit{n} n/r_n=\IF,\] then for any $\delta \in (0,\IF)$ we have
\bqn{\label{candidat}
\EIDC &:=& 
\limit{n} \frac{n}{r_n} \pk{ \max_{i\in  \{ 0, \delta \ldot \delta r_n\} } \xiW(t) > n } =
\EID \in [0,\IF).
}
\end{lemma}
In the next section we shall show that the new expression for the extremal index in \eqref{candidat}  is positive under \eqref{asC}.
  Exploiting the explicit form of the {\it tail process} we shall derive several new interesting formulas for $\HWD$. \\
Brief outline of the rest of the paper: Section 2 displays our main results which establish the positivity of the
Pickands-type constants and some new formulas.} In Section 3 we shall discuss the connection with mixed moving maxima (M3)
representation of Brown-Resnick processes. Then we derive some explicit lower bounds for $\HWD$ in case that $B$ in \eqref{Ws} is a
Gaussian or a L\'evy process and then discuss the relation between $\mathcal H_W^0$  and the
{\it mean cluster index}.
\ii{Further, we shall show that the classical Pickands constants are related to a small ball problem.} All the proofs are relegated to Section 4.

\section{Main Results}
\peH{We keep the same setup as in the Introduction and denote additionally by $\mathcal{E}$ a unit exponential random variable which is independent of everything else. According to \cite{BojanS} a candidate for the extremal index is given by the following formula
\BQN\label{formulaA}
\EIDD 
&=&  \limit{m} \pk{\max_{1 \le i \le m} Y^\delta(i)\le 1},
\EQN
where $Y^\delta(i),i\in \Z$ is the {\it tail process} of $\xiW^\delta$, see \cite{BojanS}. As in the aforementioned paper we shall impose the
{\it finite mean cluster size condition} 
of \cite[Condition 4.1]{BojanS}:
\begin{Cond} \label{C}  {\it Given $\delta>0$,   there \Ed{exists a sequence of} positive integers $r_n,n\inn$ satisfying  $\limit{n} r_n/n= \Ee{1/}\limit{n}  {r_n}=0$ such that
\BQN\label{zh}
\limit{m}\limsup_{n\to \IF} \pk{ \max_{m \le \abs{k} \le r_n} \xiW(k \delta)> nx \Ed{\Bigl \lvert \xiW(0)> nx }} =0
\EQN
holds for any $x>0$.}
\end{Cond}
 In view of \cite[Proposition 4.2]{BojanS}  \Ed{we have} that  $\EIDD>0$ \Ed{follows from} Condition \ref{C}.
Our main result below  establishes  new formulas \K{for} $\HWD$.}\\
 \ii{Moreover, from the above mentioned reference, Condition 2.1 together with well-known $\ii{\mathcal{A}}(a_n)$ conditions of
 Hsing and Davis implies that  the candidate of extremal index is equal to the extremal index, i.e., $\EIDD=\EID>0$.
 It is well-known that $\ii{\mathcal{A}}(a_n)$ is implied by the strong mixing of $\xiW^\delta$.
 However, our results derived below do not require strong mixing, but just mixing of $\xiW^\delta$.}	

\begin{theorem} \label{ThA} Let $X(t)=e^{W(t)},t\inr$ with $W$ as in \eqref{Ws} be such that  \eqref{B0} holds and
$\xiW(t),t\inr $ is max-stable and stationary.
We have that \eqref{asC}  holds for $\delta>0$ if and only if Condition \eqref{C}  holds. Moreover, if \eqref{asC} holds for $\delta>0$, then
\BQN
\HWD&=& \label{formulaAB} \peH{\frac{1}{\delta} \pk{\sup_{i< 0} \WD{i}< 0 =\sup_{i\in \Z} \WD{i}}}\\
&=& \label{albinA}  \frac{1}{\delta} \pk{\sup_{i\ge 1} \Bigl(\mathcal{E}+ \WD{i}\Bigr) \le 0}\\
&=& \label{albinB} \frac{1}{\delta} \Biggl[ \E{ \sup_{i\ge 0} e^{\WD{i}}}-\E{ \sup_{i\ge 1} e^{\WD{i}} } \Biggr]\in (0, 1/\delta),
\EQN
where $\WD{t}= W(t \delta), t\in \Z$ and $\mathcal{E}$ is a unit exponential random variable independent of $W$.
\end{theorem}
\def\HWDA{\mathcal{H}_W^0}

\begin{remark} \peH{a) If $\pk{\WD{i}=0}=0$ for any negative integer $i$,
	then
\bqny{
\pk{\sup_{- m \le i< 0} \WD{i}< 0 =\sup_{-m \le j \le m } \WD{j}} &=&
\pk{\sup_{-m \le j \le m } \WD{j}=0}
}
holds for any integer $m>1$.  Consequently, by \eqref{formulaAB} we have
\BQN
\HWD&=& \label{formulaABC}
\frac{1}{\delta} \limit{m} \pk{\sup_{- m \le i< 0} \WD{i}< 0 =\sup_{-m \le j \le m } \WD{j}} \notag\\
&=& \frac{1}{\delta} \pk{\sup_{i\in \Z} \WD{i}=0}>0,
\EQN
which has been shown  in  \cite{DiekerY} for the case $B$ is a standard fractional Brownian motion. The assumption $W(0)=0$ can be removed, see \cite{Htilt}.}\\
b) Above we assumed that $\xiW$ has c\`adl\`ag sample paths in order to define $H_W^0$. For the results of \netheo{ThA}, this assumption is not needed. \\
c) In \cite{SBK} it is shown that under the assumptions of \netheo{ThA} we have
\bqn{\label{seb}
 \HWD = \E{  \frac{  \sup_{t\in \delta \Z}  e^{ W (t)}  }  {  \delta \sum_{t\in \delta \Z} e^{W(t)} }}.
	}
According to \eqref{albinB}, for calculation of $\HWD$ it suffices to know $W(t), t\in \delta \Z, t>0$, i.e., only the values of $W$ for positive $t$ matter. This is not the case for the formula  \eqref{seb}. Both \eqref{seb}
and \eqref{albinB} are given in terms of expectations and not as  limits, which is a great advantage for simulations.
To this end, we mention that
simulation of Pickands \K{constants} has been the topic of many contributions, see e.g., \cite{Michna1,Michna2,DiekerY}.  \\
d) If $X(t)= e^{W(t)},t\inr$ is Brown-Resnick stationary, i.e., the associated max-stable process with $\zeta_W$ is max-stable and stationary, then the time reversed process $V(t)=W(-t),t\inr$ also defines a Brown-Resnick stationary processes. Moreover, for any $\delta\ge 0$
$$\HWD= \mathcal{H}_{V}^\delta.$$
 Consequently the formulas in \netheo{ThA} can be stated with $V$ instead of $W$,  for instance we have
 \bqn{
 \HWD&=& \frac{1}{\delta} \pk{\sup_{i\le -1} \Bigl(\mathcal{E}+ \WD{i}\Bigr) \le 0}\notag\\
 &=& \peH{\frac{1}{\delta} \pk{\WD{i}< 0, i\in \mathbb{N}, \WD{i} \le 0, i\in \Z }}\label{kabWang}.
}
e) If $W(t)= \sqrt{2} tL  - t^2$ with $L$ an $N(0,1)$ random variable with distribution $\Phi$ and probability density function $\varphi$, by \eqref{albinA}
we have
\bqn{  \HWD&=&   \frac{1}{\delta} \int_0^\IF \pk{\mathcal{E} + \sup_{i\ge 1} \Bigl( \sqrt{2} \delta i b  - (\delta i)^2 \Bigr)\le 0}\varphi(b) db \notag\\
	&=& \frac{1}{\delta}  \int_{- \delta/\sqrt{2}}^{ \delta/\sqrt{2}}  \varphi(b) db
	= \frac{1}{\delta}\Bigl[\Phi(\delta/\sqrt{2} )- \Phi(\K{-}\delta/\sqrt{2} )\Bigr]
}
holds for any $\delta>0$. Consequently, letting $\delta \to 0$ we obtain the well-known result
$$ \mathcal{H}_W^0= \sqrt{2} \varphi(0)=\frac{1}{\sqrt{\pi}}.$$
\end{remark}
A canonical example for $W$ with representation \eqref{Ws} is \Ee{the case when} $B$ is a centered Gaussian process with stationary increments, continuous sample paths, and  variance function $\sigma^2$. Then the max-stable process $\xiW$ is stationary,  see \cite{KabExt}.
Using a direct argument, we establish in the next theorem the positivity of $\mathcal{H}_W^0$.
 \begin{theorem} \label{c2} If
\BQN\label{ln8}
\Ee{\liminf_{t\to \IF} \frac{\sigma^2(t)}{\ln t}> 8,}
\EQN
then  $\mathcal{H}_W^0>0$.
\end{theorem}

Since \eqref{ln8} implies  \eqref{asC}, \Ee{see  Corollary 2.4 in \cite{marcus1972} or \cite{kab2009}}, then using  $\mathcal{H}_W^0\ge \HWD$ for any $\delta>0$ we immediately establish the positivity of $\mathcal{H}_W^0$.\\
Indeed, the positivity of $\mathcal{H}_W^0$ is crucial for the study of extremes of Gaussian processes.
Condition \eqref{ln8} can be easily checked, for instance if $W(t)=\sqrt{2} B_\alpha(t)- \abs{t}^\alpha$. Consequently, the classical Pickands constants
$\HWD$ are positive for any $\delta\ge 0$. This fact is highly non-trivial; after announced in Pickands' pioneering work \cite{PickandsB}, correct proofs were obtained later by Pickands himself, and in \cite{BickelRA,QuallsW}, see for instance Theorem B3  in \cite{BickelRB}.
We note in passing that  under general conditions on $\sigma^2$ the positivity of $\HWDA$ is established in
\cite{debicki2002ruin}.\\
Apart from the alternative proof for the positiveness of the original Pickands constants, \netheo{ThA} \Ee{extends to} non-Gaussian processes $W$.
For the above Gaussian setup,  direct calculations show the positivity of $\HWD$ under a slightly weaker condition than \eqref{ln8}.\\

 \def\OW{\overline{H}_W}

\section{Discussions \& Extensions}

\subsection{Relation with lower tail probabilities} For the classical case of
Piterbarg constants \KK{$\mathcal{H}_{B_\alpha}$}, i.e.,
\KK{for}
$W(t)=\sqrt{2}B_\alpha(t)- \abs{t}^\alpha, t\inr, \alpha \in (0,2]$ we show below that \eqref{formulaABC} implies
a nice relation with a small ball problem.
\begin{proposition} \label{Krzys2}
For any \KK{$\alpha\in(0,2]$} we have
\[
\lim_{\eta\to 0} \eta^{-2/\alpha} \pk{ \forall_{k\in \Z\setminus \{0\}}B_\alpha(1/k)\le\eta}
=2^{\KK{1}/\alpha} \mathcal{H}_{B_\alpha}.
\]
\end{proposition}
The above \KK{result} strongly relates to
the self-similarity property of fractional Brownian motion.
In case of \KK{a} general Gaussian $W$, we \KK{still have that} $\xiW$ is stationary
if $W$ has stationary increments.
However, fBm is the only centered Gaussian process with stationary
increments being further self-similar.
Hence, no obvious extensions of the \KK{above
relation} with lower tails can be derived for general $W$.

\subsection{Non-Gaussian $W$} The classical Pickands constants are defined for $W(t)= \sqrt{2} B_\alpha(t)- \abs{t}^\alpha$
with $B_\alpha$ a standard fBm with Hurst index $\alpha/2\in (0,1]$. The more general case where  $B_\alpha$ is substituted by a centered Gaussian process with stationary increments is discussed in details in \cite{debicki2002ruin}.\\
Our setup clearly allows for any random process $W$, not necessarily Gaussian, which is Brown-Resnick stationary.
Along with the Gaussian case of $W$, the L\'evy one has also been dealt already in the literature. In view of \cite{eng2014d, sto2008}, if $B(t), t\geq 0$ is  a L\'evy process  such that
\bqny{\label{Levi}
	\Phi(1)< \IF, \quad \Phi(\theta):=\ln \E{e^{\theta B(1)}},
}
then
$W(t)= B(t) - \Phi(1)t$, $t\geq 0$ is Brown-Resnick stationary, i.e., $\xiW(t),t\ge 0$ is max-stable stationary with unit Gumbel marginals. \\
In \cite{KW} an important constant appears in the \K{asymptotic analysis} of the maximum of standardised increments of random walks, which in fact is the Pickands constant $\HWD, \delta>0$ introduced here for $W$ as above.
In  \cite{KW}[Lemma 5.16] a new formula for $\HWD$ is derived, which is identical with our formula in \eqref{kabWang}.
Another instance of the Pickands constant given by formula \eqref{formulaAB} is displayed in \cite{SegersErgodic}[Theorem 5.3].
With the notation of that theorem, we have for $\delta=1$ that
$$W(i)= \sum_{j=1}^i A_i,$$ where $A_i$'s are iid with the same distribution as
$Z \mathbb{I}(U \le e^{- \eta Z} )$ for some $\eta>0$ with $U$ uniformly distributed on $(0,1)$ being independent of $Z$ which has some pdf symmetric around 0. \\

Pickands constants appear also in the context of semi-\K{min}-stable processes, see \cite{WangOrstein}. In view of the aforementioned paper,
 several results derived here for max-stable processes are extendable to semi-\K{min}-stable processes.

\subsection{Finite Mean Cluster Size Condition} As noted in \cite{MR2186477}, Condition \ref{C}  is implied by the
so-called {\it short-lasting exceedance condition} given below:
\begin{Cond} \label{C2}  {\it Given $\delta>0$,   there exists a sequence of integers $r_n,n\inn$ satisfying $\limit{n} r_n/n= \Ee{1/}\limit{n}{r_n}=0$ such that
\BQN\label{zhb}
\limit{m}\limsup_{n\to \IF} \sum_{k=m}^{r_n} \pk{\xiW(k \delta)> nx \Bigl \lvert \xiW(0)> nx } =0
\EQN
is valid for any $x>0$.}
\end{Cond}
This latter condition is a rephrasing of the so-called  B condition, see e.g., \cite{Albin1990,DebickiHJL14,Albin2016},
which was formulated by  discretising the original Berman's condition, see \cite{Berman82}.  Condition \ref{C2} is weaker than the $D'(x n)$ condition of Leadbetter as discussed in  \cite{Embetal1997}[Section 5.3.2].\\
Commonly, Condition \eqref{C} assumed for $x=1$ is referred to as the {\it anti-clustering condition}, see e.g.,
\cite{MR2248577,Segers16}. Clearly, the finite mean cluster size condition is stronger \K{then} the anti-clustering condition. The latter  appears in various contexts related to extremes of stationary processes,
see e.g., \cite{MR1894253,MR1805791,MR2248577,  BojanS, Segers16} and the references therein. \\

\subsection{M3 Representation}
Since we assume that $\xiW$ is max-stable stationary with c\`adl\`ag sample paths and $W$ with representation \eqref{Ws} is such that $B$ satisfies \eqref{B0},
then assuming the following almost sure convergence
\bqn {\label{fino2} W(t) \to - \IF }
as $\abs{t}\to \IF$ is equivalent with \K{the fact that} $\xiW$ possesses a {\it mixed moving maxima representation} (for short M3), \Ed{see \cite[Theorem 3]{kabDombry2} and \cite{WangStoev}}.  More specifically, under \eqref{fino2} we have the equality of finite dimensional distributions
\begin{align}\label{M3}
   \xiW (t) \stackrel{d}{=} \max_{i\ge 1} P_i e^{ F_i(t- T_i)}, \quad t\inr
\end{align}
\Ee{between rhs and lhs in \eqref{M3}}, where the $F_i$'s are independent copies of a measurable c\`adl\`ag process $F_W(t), t\inr$ satisfying
\BQN\label{normaz}
\sup_{t\in \R}  F_W(t) = F_W(0) = 0
\EQN
 almost surely, and  $\sum_{i=1}^\IF \ve_{(P_i,T_i)}$ is a PPP in $(0,\IF) \times \R$ with intensity $C_W \cdot  p^{-2}dp \cdot dt$ with
\BQN \label{cw}
C_W= \left(\E{ \int_{\R} e^{  F_W(t)} \, dt}\right)^{-1} \in (0,\IF).
\EQN
 Moreover
$\xiW^\delta$, the restriction of $\xiW$ on $\delta \Z$ possesses an M3 representation for any $\delta >0$, see \cite{SBK} for more details.  Denote the corresponding constant in the intensity of this PPP by  $C_W^\delta>0$ (and thus $C_W^0$ is just $C_W$ given in \eqref{cw}). \\
In view of \cite{SBK}[Proposition 1], if $\xiW^\delta, \delta>0$ admits an M3 representation as mentioned above, then
\bqn{ \label{fino3}
	\HWD= C_W^\delta,}
provided that \eqref{asC} holds. Hence \netheo{ThA} presents new formulas for $C_W^\delta$. \Ed{Note in passing that
\eqref{fino3} has been shown in \cite{KabExt}.  Therein it is proved that $C_W^\delta$ is given by the right-hand side of \eqref{formulaABC} assuming further that $W(t)=B(t)- \E{e^{\ln B(t)}}, t\inr$ with $B$ a centered Gaussian process with statioanry increments satisfying  $W(0)=0$ almost surely.} \\
In view of  \cite{SBK}[Theorem 1], if \eqref{asC} holds, then we have
\bqn{ \HWD= \E{ \frac{M^\delta}{S^\delta}}=C_W^\delta,}
with $M^\delta:=\max_{i\in  \Z} e^{W(i \delta)}$ and $ S^\delta:=\delta \sum_{t\in \delta \Z} e^{W(t)} $. Thus
$\HWD >0$. \\
The representation of $\HWD$ as an expectation of the ratio $M^\delta/S^\delta$  is crucial for its simulation. Such
a representation has been initially shown in \cite{DiekerY} for \K{classical Pickands constants}.

\subsection{Lower Bounds}
In \netheo{ThA} we present new formulas for $\HWD$, which in turn establish the positivity of $\HWD$ and thus the positivity for the extremal index of  $\xiW^\delta$. If only the positivity of $\HWD$ is of primary interest, then the conditions of \netheo{ThA} can be relaxed. 
Next, we consider two important classes \K{of processes} for $B$ \K{that is centered Gaussian processes with stationary increments
and L\'evy processes}.  \Ed{Results for the L\'evy case has been already given in \cite{SBK}}.\\
For particular values of $\delta$, we show that  it is possible to
derive a positive lower bound for $\HWD$ and thus establishing the positivity of $\HWD$.
\K{Let  $x_+:=\max(x,0)$.}

\begin{theorem} \label{pKRZS}
i) \K{Let $W(t)= B(t) - \sigma^2(t)/2$, $t\geq 0$, where $B(t)$ is a centered Gaussian processes with stationary increments and variance function $\sigma^2$
such that $\sigma(0)=0$. Then} for any $\delta >0$
\BQN\label{boundGA}
\HWD &\ge &
\frac{1}{\delta}  \max\Bigl(0, 1 - \sum_{k=1}^\infty e^{-\frac{\sigma^2(\delta k)}{\EHC{8}}}\Bigr).
\EQN
ii) \K{Let $W(t)= B(t) - \Phi(1)t$, $t\geq 0$,
where
$B(t)$  is a L\'evy process satisfying \eqref{Levi}. Then for any $\delta >0$}
 \BQN\label{hwLevi}
\HWD & \ge & \frac{1}{\delta}
\frac{\max\Bigl(0,1-2e^{   \left( \Phi(1/2) - \frac{1}{2} \Phi(1)\right)\delta }\Bigr)}{
1-e^{  \left( \Phi(1/2) - \frac{1}{2} \Phi(1)\right)\delta  }}.
\EQN
\end{theorem}

\begin{remark}
{
a) It follows  from i) of \netheo{pKRZS}
that if $\sigma(\delta k) \ge C (\delta k)^{\kappa/2}$ for
all $k\inn$ and some $\kappa>0$, then
\BQN
\HWD \ge  \frac{1}{\delta}  \left( 1 - \frac{1}{\delta} \frac{\Gamma(1/\kappa)}{\kappa  \left( C^2 / \EHC{8} \right)^{1/\kappa}} \right).
\EQN
}
Since $\mathcal{H}_W^0 \ge \HWD$ for any $\delta>0$, then the above implies $\mathcal{H}_W^0> 0.$\\ 
b) If $B$ is a L\'evy process as in \netheo{pKRZS}, ii), then  (see the proof in Section 3)
\BQN\label{LevyH}
 \mathcal{H}_W^0 \ge \frac{1}{8}[ \Phi(1)- 2\Phi(1/2)]> 0.
 \EQN
\end{remark}

\subsection{Case $\delta=0$}
Since \eqref{FreLim} holds also for $\delta=0$ and $\HWDA\ge \HWD$,
 then the extremal index of the continuous process $\xiW$ is
$$ \widetilde \theta _W=  \HWDA\ge 0,$$
which is positive, provided that \eqref{asC} holds. In the special case that $W(t)=\sqrt{2}B_\alpha(t)- \abs{t}^\alpha$  we have that
\BQN \label{defP2}
\lim_{\delta \downarrow 0} \HWD =\HWDA=:\HWDAA,
  \EQN
hence for such $W$ \K{and} for any $\alpha \in (0,2]$
\BQN \label{limitIndex}
 \widetilde \theta_W =\lim_{\delta \downarrow 0 } \frac{\EID}{\delta}.
 \EQN
 Recall that we denote by $\EID, \delta>0$ the extremal index of $\xiW^\delta$.
Using the terminology of \cite{H1999} we refer to $\OW$ defined by (assuming that the limit exists)
$$ \lim_{\delta \downarrow 0} \frac{\EID}{\delta}= \lim_{\delta \downarrow 0} \HWD= \OW$$
as the {\it mean cluster index of the process} $W$. Since for any $T>0$ and $\delta>0$
$$
0\le \E{\sup_{t\in  {\delta \Z} \cap [0,T]} e^{W(t)}}=:\mathcal{H}^0_W([0,T]) ,$$
 then clearly $\OW \in [0, \HWDAA]$.\\
We show next that if $\xiW$ possesses an M3 representation, then  $\OW$ is positive.
\begin{proposition}\label{prop2} Suppose that $\xiW$ is max-stable and stationary with $W(0)=0$. If $\xiW$ possesses an M3 representation and $\OW$ exists, then
\BQN \label{ow}
 \OW  \ge \E{ \frac{\sup_{t\inr} e^{W(t)}}{ \eta \sum_{t\in \eta \Z} e^{W(t)}}} > 0
 \EQN
 holds for any $\eta>0$.
\end{proposition}

\begin{remark} a) In view of \K{Theorems 2 and 3 in} \cite{SBK} we have for some general $W$ as in \eqref{Ws}, with $B$ being Gaussian or L\'evy process
\BQN \HWDA=\E{ \frac{\sup_{t\inr} e^{W(t)}}{ \eta \sum_{t\in \eta \Z} e^{W(t)}}} = \E{ \frac{\sup_{t\inr} e^{W(t)}}{ \int_{t\inr } e^{W(t)} \, dt}}
 \EQN
is valid for any $\eta>0$.  Consequently, under these conditions and the setup of \neprop{prop2}
 \BQN
\HWDA= \OW.
 \EQN
b) If $W(t)= \sqrt{2} B_\alpha(t)- \abs{t}^\alpha,t\inr $,  by \eqref{defP2} and \eqref{albinA} for any
$\alpha \in (0,2]$
\BQN \label{albinAA}
\HWDA &=& \OW
= \lim_{\delta \downarrow 0}  \frac{1}{\delta} \pk{\sup_{i\ge 1} (\mathcal{E}+ \WD{i}) \le 0} ,
\EQN
with $\mathcal{E}$ a unit exponential random variable independent of $W$.
\K{Expression \eqref{albinAA} of the classical Pickands constant was} initially derived in \cite{Albin1990} for some general $W$,
see also \K{recent contribution} \cite{Albin2016}.
In \cite{H1999}, Proposition 3 or the formula in \cite{FalkHusler}[p.44] the classical Pickands constant  is the  limit of a cluster index.\\
\COM{c) For $B$ a fractional Brownian motion \cite{DiekerY} proved various formulas for $\HWD$ and showed that these constants can be expressed in terms of expectations. Our new expression of $\HWD$ in \eqref{albinB} is also given in terms of expectations, however it seems to be unrelated to the representations derived in the aforementioned paper.
}
\end{remark}

\COM{
\subsection{Comparison criterion for extremal indices and Slepian-type inequalities}

In this section we focus on the case where
$W(t)=B(t) -\sigma^2(t)/2$, with $B(t)$ a sample continuous centered Gaussian
process with stationary increments and
variance function $\sigma^2(t)$  such that
\\
\\
{\bf{A1}} $\sigma_i^2(t)\le |t|^\alpha$ as $t\to 0$  for some $\alpha>0$.
\\
\\
The following finding provides an extension of Theorem 3.2 in \cite{debicki2002ruin}.
\begin{theorem}\label{th.pick}
Suppose that $W_i(t),$ $i=1,2$ satisfy {\bf A1}.
If for a given compact set
$E\subset \R$ and all $t\in E$
\[
\sigma_1^2(t)\ge\sigma_2^2(t),
\]
then
\[
\E{ \sup_{t\in E}e^{ W_1(t) }}\ge \E{\sup_{t\in E}e^{ W_2(t) }}.
\]
\end{theorem}
\begin{remark}
From the proof of, e.g., Theorem 2.1 in \cite{DHL16} it follows that in the case of $E=\delta \Z \cap G$ with $\delta>0$ and compact $G\in R$,
Theorem \ref{th.pick} holds also without assumption ${\bf A1}$.
\end{remark}

Theorem \ref{th.pick} combined with
the fact
that for any compact $E\in\R$ and $x\in\R$
\[
\pk{ \xiW(t) \le x, \forall t\in E} = e^{- \E{\max_{t\in E} \bigl(e^{W(t)}/x \bigr)}}
\]
allows
us to obtain both comparison criterion for
extremal indices  of $\xi_W^\delta(t)$, $\delta>0$, and Slepian-type inequalities
for  $\xi_W(t)$.

\begin{proposition}
Suppose that, for $\delta>0$,
\[
\sigma_1^2(\delta i)\ge\sigma_2^2(\delta i),
\]
for all $i\in \N$.
Then extremal indices $\theta^\delta_{W_i}$ of max-stable processes $\xi^\delta_{W_i}$, $i=1,2$ satisfy
\[
\theta^\delta_{W_1}\ge \theta^\delta_{W_2}.
\]
\end{proposition}

\begin{proposition}
Suppose that $W_i(t),\ t\in\R,$ $i=1,2$ satisfy {\bf A1}.
If for a given compact set
$E\subset \R$ and all $t\in E$
\[
\sigma_1^2(t)\ge\sigma_2^2(t),
\]
then for each $u\in \R$
\[
\pk{ \sup_{t\in E} \xi_{W_1}(t)>u} \ge \pk{ \sup_{t\in E} \xi_{W_2}(t)>u}.
\]
\end{proposition}
}

\section{Proofs}

\peH{\prooflem{propAA} Since $\limit{n} r_n= \IF$, then by \eqref{FreLim} and \eqref{FreLim2} 
$$
\limit{n} r_n^{-1}\E{ \max_{i\in  \{ 0, \delta \ldot \delta r_n\} }  e^{W(i)}}= \delta \HWD= \EID.
$$
For any $n\inn$ we have
\bqny{
\lefteqn{\frac{ \pk{ \max_{i\in  \{ 0, \delta \ldot \delta r_n\} } \xiW(i) > n }}{r_n \pk{\xiW(0)> n}}}\\
&=& \frac{  \pk{ \max_{i\in  \{ 0, \delta \ldot \delta r_n\} } \xiW(i) >  n }}{r_n [1- e^{-1/n}]}\\
&\sim& n r_n^{-1} \Bigl[ 1 - \pk{ \max_{i\in  \{ 0, \delta \ldot \delta r_n\} } \xiW(i) \le  n } \Bigr]\\
&=& n r_n^{-1}\Bigl[1 - e^{ - c_n/n} \Bigr], \quad c_n:= \E{ \max_{i\in  \{ 0, \delta \ldot \delta r_n\} }  e^{W(i)}},
}
where the last equality follows from \eqref{base}.
The assumption that $\limit{n} n/r_n=\IF$ and $\E{e^{W(i)}}=1, i\in \delta \Z$ imply
\bqn{ \label{cnD}
	\frac{c_n}n \le \frac 1 n \E{ \sum_{ i\in  \{ 0, \delta \ldot \delta r_n\} } e^{W(i)} } = \frac{r_n+1}{n} \to 0, \quad n\to \IF.
	}
Consequently,
\bqny{
\frac{ \pk{ \max_{i\in  \{ 0, \delta \ldot \delta r_n\} } \xiW(i) > n }}{r_n \pk{\xiW(0)> n}}
	&\sim &  r_n^{-1}\E{ \max_{i\in  \{ 0, \delta \ldot \delta r_n\} }  e^{W(i)}}
	\sim  \EID, \quad n\to \IF,
}
hence the claim follows.  \QED
}

\prooftheo{ThA} We show first stochastic representation \eqref{stochRep}. Recall that $X(t)= e^{W(t)}$ and
for $\delta>0$ we set
$$\WD{t}= W(\delta t), \quad
\XD{t}=e^{\WD{t}}, \quad t\in \Z.$$
 By \eqref{base},
the fact that $\pk{\xiW(0)\le x}=e^{-1/x}, x>0$  and the assumption that $X(0)=1$  almost surely, for any $y_1 \ldot y_n$ positive and $y_0 > 1$ we have
{\tiny
\BQNY
\lefteqn{ \pk{ \xiWD{i}  \le T  y_i, i=0 \ldot n \Bigl \lvert \xiWD{0}> T  }}\\
&=&
\frac{ 1- \pk{\xiWD{0} \le T, \xiWD{i} \le T y_i, i \in \{0 \ldot n\} } - [1-
 \pk{ \xiWD{i} \le T y_i,i \in \{0 \ldot n\} }]  }{\pk{ \xiWD{0}> T}}
 \\
&=& \frac{ 1- e^{- \EE{ \max\Bigl( \XD{0},    \max_{i \in \{1 \ldot n\}} \frac{\XD{i}}{y_i}\Bigr)} \frac{1}T }  -
 \Biggl[1- e^{- \EE{ \max_{i \in \{1 \ldot n\}} \frac{\XD{i}}{y_i}} \frac{1}T }\Biggr]  }{1- e^{-\frac 1 T} }\\
&\sim & T \Biggl[ 1-  \biggl[ 1- \frac{1}T\mathbb{E}\biggl\{ \max \Bigl( 1,    \max_{i \in \{0 \ldot n\}} \frac{\XD{i}}{y_i}\Bigr) \biggr\}\biggl]
- \Biggl(1-   \biggl[ 1- \frac{1}T\EE{ \max_{i \in \{0 \ldot n\}} \frac{\XD{i}}{y_i}}\biggl]   \Biggr) \Biggr]\\
&\to & \E{\Biggl(1-   \max_{i \in \{ 0\ldot n\}} \frac{\XD{i}}{y_i}\Biggr)_+}, \quad T \to \IF\\
&=& \pk{ \ii{\mathcal{P}}   \le y_0, \ii{\mathcal{P}}  \XD{i} \le y_i,  \forall i\in \{1 \ldot n\}},
\EQNY
}
where $\ii{\mathcal{P}} $ is a unit Pareto random variable with survival function $1/s, s>1$ independent of the process $X$.  Hence the claim in \eqref{stochRep} follows by \cite{BojanS}[Theorem 2.1 (ii)]. Next 
by the above derivations for any 
sequence of integers $r_n> m \inn$ for any $x>0$ (recall $X^\delta(0)=1$ almost surely) we have
{\tiny
\BQNY
\lefteqn{1- \pk{ \max_{  m \le \abs{i} \le r_n} \xiWD{i} > nx  \Bigl \lvert \xiWD{0} > nx }}\\
&=&
\frac{\pk{ \max_{  m \le \abs{i}  \le r_n} \xiWD{i} \le  nx ,  \xiWD{0} > nx }}{\pk{\xiWD{0} > nx}}\\
&=&  \frac{ 1- e^{- \E{ \max( \XD{0},    \max_{\abs{i} \in \{m \ldot r_n\}} \XD{i})} \frac{1}{nx} }  -
	\Bigl[1- e^{- \E{ \max_{\abs{i} \in \{m \ldot r_n\}} \XD{i} } \frac{1}{nx} }\Bigr]  }{1- e^{-\frac 1 {nx}} }\\
&\sim & nx \Biggl[ 1-  \Bigl[ 1- \frac 1 {nx}\E{ \max \Bigl(  1,    \max_{\abs{i} \in \{m \ldot r_n\}} \XD{i} \Bigr)}\Bigl]
- \Biggl(1-   \Bigl[ 1- \frac 1 {nx}\E{ \max_{\abs{i} \in \{m \ldot r_n\}}  \XD{i} }\Bigl]   \Biggr) \Biggr]\\
&\sim& \E{\Bigl(1-   \max_{\abs{i} \in \{ m\ldot r_n\}} \XD{i}\Bigr)_+},
\EQNY
}
where we used the fact that as in \eqref{cnD}, the condition  $\limit{n}r_n= \limit{n}\frac{n}{r_n}= \IF$  implies
$$ \limit{n} \frac{1} n \E{ \max\Bigl( {\XD{0}} ,    \max_{\abs{i} \in \{m \ldot r_n\}} {\XD{i}} \Bigr)} =0,
$$
and
$$ \limit{n} \frac{1}n \E{ \max_{\abs{i} \in \{m \ldot r_n\}} {\XD{i}}}  =0.$$
Consequently,
 \BQNY
\lefteqn{\limit{m} \limsup_{n\to \IF} \pk{ \max_{  m \le \abs{i} \le r_n} \xiWD{i} > nx  \Bigl \lvert \xiWD{0} > nx }=}\nonumber\\
&=& \limit{m} \limsup_{n\to \IF}   \Biggl[ 1-
\E{\Bigl(1-   \max_{\abs{i} \in \{ m\ldot r_n\}} \XD{i}\Bigr)_+}\Biggr]\\
&=& 1- \limit{m}
\E{\Bigl(1-   \max_{\abs{i} \in \Z, i\ge m} \XD{i}\Bigr)_+}\\
&=&0,
\EQNY
\peH{where we used the assumption \eqref{asC}. Hence Condition \ref{C} holds.\\
In light of \cite[Proposition 4.2]{BojanS} we have that Condition \ref{C} implies \eqref{asC}.  Moreover, since
$$ \pk{\xiW(0)> n}= 1- e^{-1/n} \sim \frac 1 n, \quad n\to \IF$$
\cite[Proposition 4.2]{BojanS}  and \nelem{propAA} imply
\BQNY
\EID=\EIDC=\EIDD > 0.
\EQNY
}
Consequently,
\BQN
 \EIDD&=& \pk{ \sup_{i\ge 1} \YD{i} \le 1} \notag\\
 & =&  \limit{n} \pk{ \K{\ii{\mathcal{P}}} \sup_{n \ge i \ge 1} \XD{i} \le 1} \label{thA1}\\
& = & \limit{n} \E{ \Bigl(1- \sup_{ n\ge i\ge 1} \XD{i} \Bigr)_+} \notag\\
& = & \E{ \Bigl (1- \sup_{ i\ge 1} \XD{i} \Bigr)_+} \notag\\
&=&  \E{ \sup_{i\ge 0} \XD{i}- \sup_{i\ge 1} \XD{i} } \in (0,1], \notag
\EQN
where the second last expression follows from the monotone convergence theorem. In fact, the above claim readily follows also from \cite{BojanS}[Remark 4.7].
Further from \eqref{thA1} we obtain
\BQNY
 \limit{n} \pk{ \ii{\mathcal{P}} \sup_{n \ge i \ge 1} \XD{i} \le 1} &=&
\limit{n} \pk{  \sup_{n \ge i \ge 1}\Bigl ( \ln \ii{\mathcal{P}} + \ln \XD{i} \Bigr) \le 0}\\
&=& \limit{n} \pk{  \sup_{n \ge i \ge 1} \Bigl( \mathcal{E}+ \WD{i} \Bigr) \le 0} \\
& =&   \pk{  \sup_{  i \ge 1} \Bigl( \mathcal{E}+ \WD{i} \Bigr) \le 0},
\EQNY
\COM{ \peH{In view of \cite[Theorem 4.3]{BojanS} we have
$$ \EIDC= \pk{ \sup_{i< - 1} e^{ W^\delta(i)} \le 1}=
\pk{ \sup_{i< - 1}  W^\delta(i) \le 0}$$
and by \cite[Theorem 2.5]{BojanSArxiv}
$$ \EIDC =\pk{ \sup_{i< - 1} e^{ W^\delta(i)} <  e^{W^\delta(0)}=\sup_{i \in \Z }  e^{ W^\delta(i)} }  $$
or equivalently since $W(0)=0$
$$ \EIDC =\pk{ \sup_{i< - 1} W^\delta(i)  <  0 = \sup_{i \in \Z } W^\delta(i) }  $$
}
}
with $\mathcal E= \ln \mathcal{P}$ a unit exponetial random variable independent of $X$. \\
Next, \eqref{formulaAB} follows from \cite{MR2186477}[Eq. (16)]. Since further we assume   \eqref{Ws}, then
\eqref{formulaAB}   implies
\bqn{ \HWD \in (0, 1/\delta)}
for any $\delta>0$,  establishing thus the proof. \QED

\prooftheo{c2} By our assumption for all large $k$
$$ \frac{\sigma^2(\delta k)}{8} > \ln (\delta k)^a. $$
Consequently, by \eqref{boundGA} we have for all $\delta$ large and some $a>1$
\BQNY
\mathcal{H}_W^0 \ge \HWD  &\ge & \frac{1}{\delta}  \left( 1 - \sum_{k=1}^\infty e^{-
	\frac{ \sigma^2(\delta k) }{\EHC{8}}}\right)
\ge
\frac{1}{\delta}  \left( 1 - \frac{1}{\delta^a}\sum_{k=1}^\infty  \frac{1}{k^a}\right)>0.
\EQNY
Hence the proof is complete. \QED

\proofprop{Krzys2}
Since $B_\alpha(0)=0$ almost surely, in view of \KK{\eqref{formulaABC}}
(see also \cite{DiekerY}[Proposition 4]) we obtain
\[
\lim_{\delta\to 0} \delta^{-1}
{\pk{ \forall_{k \in \Z\setminus \{0\}}B_\alpha(\delta k)\le |\delta k|^\alpha/\KK{\sqrt{2}}}}
=\mathcal{H}_{B_\alpha}.
\]
Moreover, by the self-similarity of $B_\alpha$, we have
\begin{eqnarray*}
	\pk{ \forall_{k\in \Z\setminus \{0\}}B_\alpha(\delta k)\le |\delta k|^\alpha/\KK{\sqrt{2}}}
	&=&
	\pk{ \forall_{k\in \Z\setminus \{0\}}|\delta k|^\alpha B_\alpha\left(\frac{1}{\delta k}\right)
\le |\delta k|^\alpha/\KK{\sqrt{2}}}\\
	&=&
	\pk{ \forall_{k\in \Z\setminus \{0\}}B_\alpha\left(\frac{1}{\delta k}\right)\le 1/\KK{\sqrt{2}}}\\
	&=&
	\pk{ \forall_{k\in \Z\setminus \{0\}}B_\alpha\left(\frac{1}{ k}\right)\le \delta^{\alpha/2}/\KK{\sqrt{2}}},
\end{eqnarray*}
hence the proof follows easily. \QED

\prooftheo{pKRZS}
i) The proof is based on  a technique
developed in Lemma 16 and  Corollary 17 in \cite{demiro2003simulation} and in Lemma 7 in \cite{shao1996bounds}, therefore we omit some details.
For any $\delta>0$ and $T$ positive integer, using Bonferoni's inequality we have  for any
process $W$ such that $\E{e^{W(k\delta)}}=1, k\ge 1$ 
\BQN
\lefteqn{\E{ \sup_{t\in \delta \Z \cap [0,\delta T]} e^W(t)}}\notag\\
&=& \int_{\R}  e^s \pk{\sup_{t\in \delta \Z \cap [0,\delta T]} W(t)> s  }\, ds\notag\\
&\ge &  \int_{\R}  e^s \pk{\exists_{1 \le  k \le T} W(k\delta)> s} \, ds \notag\\
&\ge &  \sum_{k=1}^{T}\int_{\R}  e^s \pk{ W(k\delta)> s} \, ds
\nonumber\\
&& -\ \  \sum_{k=1}^{T-1} \sum_{l=k+1}^T \int_{\R} e^s \pk{ W(k\delta)> s, W(l \delta)> s } \, ds \notag\\
&\ge &  \sum_{k=1}^{T}\E{ e^{W(k\delta)}} \sum_{k=1}^{T-1} \sum_{l=k+1}^T \int_{\R} e^s \pk{ W(k\delta)+  W(l \delta)> 2s } \, ds \notag\\
&=&  T- \sum_{k=1}^{T-1} \sum_{l=k+1}^T \int_{\R} e^s \pk{ W(k\delta)+  W(l \delta)> 2s } \, ds \label{30}\\
&= &T -\sum_{k=1}^{T-1} \sum_{l=k+1}^T \E{e^{\frac{W(k\delta)+ W(l\delta)}{2}}}  \notag\\
&= &T -\sum_{k=1}^{T-1} \sum_{l=k+1}^T e^{- \frac{\sigma^2(\delta \abs{k-l})}{\EHC{8}}}\notag\\
&\ge &
\KD{T-T\sum_{k=1}^T e^{- \frac{\sigma^2(\delta k)}{\EHC{8}}}}, \notag
%
\EQN
where the last equality follows by the stationary of increments of the random process $B$. Along the lines of the proof in \cite{Tabis}
\KD{
\begin{eqnarray*}
\mathcal{H}_W^\delta&=& \limit{T} \frac{1}{T}\E{ \sup_{t\in \delta \Z \cap [0,T]} e^W(t)}\\
&\ge & \limit{T}  \frac{1}{T}  \floor{T/\delta}\Bigl(1-
\sum_{k=1}^\IF e^{- \frac{\sigma^2(\delta k)}{\EHC{8}}}\Bigr)_+ \\
&=& \frac{1}{\delta} \Bigl(1-
\sum_{k=1}^\IF e^{- \frac{\sigma^2(\delta k)}{\EHC{8}}}\Bigr)_+.
\end{eqnarray*}
\COM{Further, if
$\sigma(t) \ge C t^{\kappa/2}$ for some $\kappa>0$, then due to the above
$\mathcal{H}_W^\delta
\ge
\frac{1}{\delta}  \left( 1 - \frac{1}{\delta} \frac{\Gamma(1/\kappa)}{\kappa  \left( C^2 / 4 \right)^{1/\kappa}} \right).
$
}
}
ii) \K{In} view of  \eqref{30},
in order to establish the proof we need to calculate
$$a_{kl}=  \int_\R e^s \pk{ W(\delta k)+W(\delta l)>2s}\, ds.$$
\COM{

For a fixed $ \delta > 0 $ and a positive sufficiently large integer  $ N $, by Bonferroni inequality,
we have
\begin{eqnarray*}
\lefteqn{
\int_\R e^w \pk{\sup_{t\in \delta \Z \cap [0,N]}( B^+(t)-\Phi(1)t)>w} dw
\ge}\\
&\ge&
\int_\R e^w \pk{\exists_{1\le k \le N/\delta}( B^+(\delta k)-\Phi(1)\delta k)>w} dw\\
&\ge&
\sum_{k=1}^{N/ \delta }
 \int_\R e^w \pk{ B^+(\delta k)-\Phi(1)\delta k>w} dw -\\
&&-
\sum_{k=1}^{N/ \delta -1}
\sum_{l=k+1}^{N/ \delta}
 \int_\R e^w \pk{ B^+(\delta k)+B^+(\delta l)-\Phi(1)\delta (k+l)>2w} dw
\end{eqnarray*}
Observe that for each $k<l$
\[
 \int_\R e^w \pk{ B^+(\delta k)-\Phi(1)\delta k>w} dw=1
\]
and further}
By the independence of the increments, and the fact that $W(\delta l)- W(\delta k) \EQD W(\delta (l-k))$ we have
\begin{eqnarray*}
a_{kl}&=&\E{e^{ \frac{W(\delta k)+ W(\delta l )}{2}}}\\
&=&\E{ e^{W(\delta k)}} \E{e^{ \frac{ W(\delta l)- W(\delta k)}{2}}}\\
&=&\E{ e^{W(\delta k)}} \E{e^{ \frac{ W(\delta (l-k))}{2}}}\\
&=&
\E{e^{\frac{B(\delta (l-k))-\Phi(1)\delta (l-k)}{2}}}\\
&=&
\exp \left( - \delta(l-k) \lambda \right),
\end{eqnarray*}
where \ii{$\lambda :=  \frac{1}{2} \Phi(1)-  \Phi(1/2)>0$
by Jensen's inequality and independence
and stationarity of increments of the L\'evy process $B$.
} Consequently, for $N\inn$ we obtain
\begin{eqnarray}
\int_\R e^s \pk{\sup_{t\in \delta \Z \cap [0,N]} W(t)>s}\, ds
&\ge&
\frac{N}{\delta}
\left( 1- \sum_{k=1}^{\infty} e^{ -\delta k \lambda } \right)\label{for.1}\\
&=&
\frac{N}{\delta}
\frac{1-2\exp \left( - \delta  \lambda  \right)}{
1-\exp \left( -\delta  \lambda  \right)},\nonumber
\end{eqnarray}
which leads to
\[
\mathcal{H}_W^\delta
\ge
\frac{1}{\delta}
\frac{1-2\exp \left( - \delta  \lambda  \right)}{
	1-\exp \left( -\delta  \lambda  \right)}
\]
and thus the proof is complete. \QED

{\bf Proof of \eqref{LevyH}}:
By (\ref{for.1}) and letting \ii{$\lambda= \frac{1}{2} \Phi(1)- \Phi(1/2)  >0$} we have
\begin{eqnarray*}
\HWDA
&\ge&
\limit{N}
\frac 1 N \int_\R e^s \pk{\sup_{t\in \delta \Z \cap [0,N]} W(t) >s}\,  ds\\
&\ge&
\frac{1}{\delta}\left(1-   \sum_{k=1}^{\infty} e^{ - \delta k \lambda  } \right) \\
&\ge&
\frac{1}{\delta}\left(1-
 \int_0^\infty
e^{ - \delta x \lambda  } dx\right)\\
&=&
\frac{1}{\delta}\left( 1-\frac{1}{\delta   \lambda }   \right)\\
&\ge& \frac \lambda  4 >0
\end{eqnarray*}
establishing the proof. \QED

\proofprop{prop2} In view of \cite{SBK} for any $\delta>0$ and any integer $k\inn$ we have
$$\HWD \ge \E{ \frac{\sup_{t\in \delta \Z} e^{W(t)}}{
k \delta \sum_{t\in k \delta \Z}  e^{W(t)}}},$$
hence choosing   $\delta_n= \eta l^{-n}$ with $\eta>0$ and $l>1$ some integer and for
$k= l^n$ which is clearly integer for any $n\ge 1$ we have
\BQNY
\mathcal{H}_W^{\delta_n} &\ge &
\E{ \frac{\sup_{t\in \delta_n \Z} e^{W(t)}}{k \delta_n \sum_{t\in k \delta_n \Z}  e^{W(t)}}}\\
&=&  \E{ \frac{\sup_{t\in \delta_n \Z} e^{W(t)}}{\eta  \sum_{t\in \eta   \Z}  e^{W(t)}}}\\
&\to & \E{ \frac{\sup_{t\inr} e^{W(t)}}{\K{\eta} \sum_{t\in \eta  \Z}  e^{W(t)}}}, \quad n\to \IF,
\EQNY
where the last limit follows by the monotone convergence theorem and the fact that $W$ has continuous sample paths.
Since by the construction $\mathcal{H}_W^{\delta_n}$ is non-decreasing in $n$, and we assume that
$\lim_{\delta \downarrow 0} \HWD=\OW$, then the claim follows.
\QED
\\

{\bf Acknowledgements.}
We are indebted to Bojan Basrak and Philippe Soulier for several comments and suggestions.
Thanks to  a Swiss National Science Foundation grant  no.  200021-166274 for partial support.
KD also acknowledges partial support by NCN Grant No 2015/17/B/ST1/01102 (2016-2019).


\begin{thebibliography}{99}                  
\def\polhk#1{\setbox0=\hbox{#1}{\ooalign{\hidewidth
  \lower1.5ex\hbox{`}\hidewidth\crcr\unhbox0}}}
  \def\polhk#1{\setbox0=\hbox{#1}{\ooalign{\hidewidth
  \lower1.5ex\hbox{`}\hidewidth\crcr\unhbox0}}}


\bibitem{Albin1990}
\n{ J.~M.~P. Albin}, {\em On extremal theory for stationary processes}, Ann.
  Probab., 18 (1990), pp.~92--128.

\bibitem{Albin2016}
\n{ J.~M.~P. Albin, E.~Hashorva, L.~Ji, and C.~Ling}, {\em Extremes and limit
  theorems for difference of chi-type processes}, ESAIM: P\& S, 20, (2016), pp.~349--366.

\bibitem{MR1894253}
\n{ B.~Basrak, R.~A. Davis, and T.~Mikosch}, {\em Regular variation of {GARCH}
  processes}, Stochastic Process. Appl., 99 (2002), pp.~95--115.

\bibitem{MR3025708}
\n{ B.~Basrak, D.~Krizmani{\'c}, and J.~Segers}, {\em A functional limit
  theorem for dependent sequences with infinite variance stable limits}, Ann.
  Probab., 40 (2012), pp.~2008--2033.

\bibitem{BojanS}
\n{ B.~Basrak and J.~Segers}, {\em Regularly varying multivariate time
  series}, Stochastic Process. Appl., 119 (2009), pp.~1055--1080.

\bibitem{Berman82}
\n{ S.~Berman}, {\em Sojourns and extremes of stationary processes}, Ann.
  Probab., 10 (1982), pp.~1--46.

\bibitem{BickelRA}
\n{ P.~J. Bickel and M.~Rosenblatt}, {\em On some global measures of the
  deviations of density function estimates}, Ann. Statist., 1 (1973),
  pp.~1071--1095.

\bibitem{BickelRB}
\leavevmode\vrule height 2pt depth -1.6pt width 23pt, {\em Two-dimensional
  random fields}, in Multivariate analysis, {III} ({P}roc. {T}hird {I}nternat.
  {S}ympos., {W}right {S}tate {U}niv., {D}ayton, {O}hio, 1972), Academic Press,
  New York, 1973, pp.~3--15.

\bibitem{Michna1}
\n{ K.~Burnecki and Z.~Michna}, {\em Simulation of {P}ickands constants},
  Probab. Math. Statist., 22 (2002), pp.~193--199.

\bibitem{MIKYuw}
\n{ R.~A. Davis, T.~Mikosch, and Y.~Zhao}, {\em Measures of serial extremal
  dependence and their estimation}, Stochastic Process. Appl., 123 (2013),
  pp.~2575--2602.

\bibitem{SBK}
\n{ K.~D\c{e}bicki, S.~Engelke, and E.~Hashorva}, {\em Brown-{R}esnick
  processes and {P}ickands-type constants}, Extremes, DOI: 10.1007/s10687-017-0289-1,  (2017).

\bibitem{DebickiHJL14}
\n{ K.~D\c{e}bicki, E.~Hashorva, L.~Ji, and C.~Ling}, {\em Extremes of order
  statistics of stationary processes}, TEST, 24 (2015), pp.~229--248.

\bibitem{debicki2002ruin}
\n{ K.~D{\c{e}}bicki}, {\em Ruin probability for {G}aussian integrated
  processes}, Stochastic Process. Appl., 98 (2002), pp.~151--174.

\bibitem{Tabis}
\n{ K.~D{\c{e}}bicki, E.~Hashorva, L.~Ji, and K.~Tabi{\'s}}, {\em Extremes of
  vector-valued {G}aussian processes: {E}xact asymptotics},
Stochastic Process.  Appl., 125 (2015), pp.~4039--4065.

\bibitem{KrzKos}
\n{ K.~D{\polhk{e}}bicki and K.~M. Kosi{\'n}ski}, {\em An
  {E}rd\"os-{R}\'ev\'esz type law of the iterated logarithm for order
  statistics of a stationary {G}aussian process}, J. Theor. Probab.,  (2016).

\bibitem{demiro2003simulation}
\n{ K.~D{\c{e}}bicki, Z.~Michna, and T.~Rolski}, {\em Simulation of the
  asymptotic constant in some fluid models}, Stochastic Models, 19 (2003),
  pp.~407--423.

\bibitem{Delrome}
\n{ M.~Delorme, A.~Rosso, and K.~J. Wiese}, {\em Pickands' constant at first
  order in an expansion around {B}rownian motion},
  https://arxiv.org/abs/1609.07909,  (2016).

\bibitem{DM}
\n{ A.~B. Dieker and T.~Mikosch}, {\em Exact simulation of {B}rown-{R}esnick
  random fields at a finite number of locations}, Extremes, 18 (2015),
  pp.~301--314.

\bibitem{DiekerY}
\n{ A.~B. Dieker and B.~Yakir}, {\em On asymptotic constants in the theory of
  extremes for {G}aussian processes}, Bernoulli, 20 (2014), pp.~1600--1619.

\bibitem{kabDombry2}
\n{ C.~Dombry and Z.~Kabluchko}, {\em Ergodic decompositions of stationary
  max-stable processes in terms of their spectral functions},
Stochastic Process.  Appl., DOI: 10.1016/j.spa.2016.10.00 (2017,in press).

\bibitem{Jakub}
\n{ P.~Doukhan, A.~Jakubowski, and G.~Lang}, {\em Phantom distribution
  functions for some stationary sequences}, Extremes, 18 (2015), pp.~697--725.

\bibitem{Embetal1997}
\n{ P.~Embrechts, C.~Kl{\"u}ppelberg, and T.~Mikosch}, {\em Modelling extremal
  events for insurance and finance}, Springer-Verlag, Berlin, 1997.

\bibitem{eng2014d}
\n{ S.~Engelke and Z.~Kabluchko}, {\em Max-stable processes associated with
  stationary systems of independent {L}{\'e}vy particles}, Stochastic Process.
  Appl., 125 (2015), pp.~4272--4299.

\bibitem{FalkHusler}
\n{ M.~Falk, J.~H{\"u}sler, and R.-D. Reiss}, {\em Laws of small numbers:
  extremes and rare events}, Birkh\"auser/Springer Basel AG, Basel,
  extended~ed., 2011.

\bibitem{Harper1}
\n{ A.~J. Harper}, {\em Bounds on the suprema of {G}aussian processes, and
  omega results for the sum of a random multiplicative function}, Ann. Appl.
  Probab., 23 (2013), pp.~584--616.

\bibitem{Harper2}
\leavevmode\vrule height 2pt depth -1.6pt width 23pt, {\em Pickands' constant
  ${H}_{\alpha}$ does not equal $1/{\Gamma}(1/\alpha)$, for small $\alpha$},
  Bernoulli, accepted,  (2015).

\bibitem{Htilt}
\n{ E.~Hashorva}, {\em Representations of max-stable processes via exponential
  tilting}, arXiv:1605.03208v3,  (2016).

\bibitem{H1999}
\n{ J.~H{\"u}sler}, {\em Extremes of a {G}aussian process and the constant
  {$H_\alpha$}}, Extremes, 2 (1999), pp.~59--70.

\bibitem{MR3301293}
\n{ A.~Janssen and J.~Segers}, {\em Markov tail chains}, J. Appl. Probab., 51
  (2014), pp.~1133--1153.

\bibitem{kab2009}
\n{ Z.~Kabluchko, M.~Schlather, and L.~de~Haan}, {\em Stationary max-stable
  fields associated to negative definite functions}, Ann. Probab., 37 (2009),
  pp.~2042--2065.

\bibitem{KW}
\n{ Z.~Kabluchko and Y.~Wang}, {\em Limiting distribution for the maximal
  standardized increment of a random walk}, Stochastic Process. Appl., 124
  (2014), pp.~2824--2867.

\bibitem{Danijel}
\n{ D.~Krizmani{\'c}}, {\em Functional weak convergence of partial maxima
  processes}, Extremes, 19 (2016), pp.~7--23.


\bibitem{Soulier}
\n{ R.~Kulik and P.~Soulier}, {\em Heavy tailed time series with extremal
  independence}, Extremes, 18 (2015), pp.~273--299.

\bibitem{marcus1972}
\n{ M.~Marcus}, {\em Upper bounds for the asymptotic maxima of continuous
  {G}aussian processes}, Ann. Math. Statist., 43 (1972), pp.~522--533.

\bibitem{MR3179972}
\n{ N.~M. Markovich}, {\em Modeling clusters of extreme values}, Extremes, 17
  (2014), pp.~97--125.

\bibitem{Michna2}
\n{ Z.~Michna}, {\em On tail probabilities and first passage times for
  fractional {B}rownian motion}, Math. Methods Oper. Res., 49 (1999),
  pp.~335--354.

\bibitem{MR1805791}
\n{ T.~Mikosch and C.~St{\u{a}}ric{\u{a}}}, {\em Limit theory for the sample
  autocorrelations and extremes of a {GARCH} {$(1,1)$} process}, Ann. Statist.,
  28 (2000), pp.~1427--1451.

\bibitem{MikoshZ}
\n{ T.~Mikosch and Y.~Zhao}, {\em The integrated periodogram of a dependent
  extremal event sequence}, Stochastic Process. Appl., 125 (2015),
  pp.~3126--3169.

\bibitem{MolchanovSPA}
\n{ I.~Molchanov and K.~Stucki}, {\em Stationarity of multivariate particle
  systems}, Stochastic Process. Appl., 123 (2013), pp.~2272--2285.

\bibitem{KabExt}
\n{ M.~Oesting, Z.~Kabluchko, and M.~Schlather}, {\em Simulation of
  {B}rown-{R}esnick processes}, Extremes, 15 (2012), pp.~89--107.

\bibitem{PickandsB}
\n{ J.~Pickands, III}, {\em Asymptotic properties of the maximum in a
  stationary {G}aussian process}, Trans. Amer. Math. Soc., 145 (1969),
  pp.~75--86.

\bibitem{Pit20}
\n{ V.~I. Piterbarg}, {\em Twenty Lectures About {G}aussian Processes},
  Atlantic Financial Press, London, New York, 2015.

\bibitem{QuallsW}
\n{ C.~Qualls and H.~Watanabe}, {\em Asymptotic properties of {G}aussian
  processes}, Ann. Math. Statist., 43 (1972), pp.~580--596.

\bibitem{SegersErgodic}
\n{ G.~O. Roberts, J.~S. Rosenthal, J.~Segers, and B.~Sousa}, {\em Extremal
  indices, geometric ergodicity of {M}arkov chains, and {MCMC}}, Extremes, 9
  (2006), pp.~213--229.

\bibitem{MR2186477}
\n{ J.~Segers}, {\em Approximate distributions of clusters of extremes},
  Statist. Probab. Lett., 74 (2005), pp.~330--336.

\bibitem{MR2248577}
\leavevmode\vrule height 2pt depth -1.6pt width 23pt, {\em Rare events,
  temporal dependence, and the extremal index}, J. Appl. Probab., 43 (2006),
  pp.~463--485.

\bibitem{Segers16}
\n{ J.~Segers, Z.~Yuwei, and M.~Thomas}, {\em Radial-angular decomposition of
  regularly varying time series in star-shaped metric spaces},
  https://arxiv.org/abs/1604.00241,  (2016).

\bibitem{shao1996bounds}
\n{ Q.~Shao}, {\em Bounds and estimators of a basic constant in extreme value
  theory of {G}aussian processes}, Statistica Sinica, 6 (1996), pp.~245--258.

\bibitem{sto2008}
\n{ S.~A. Stoev}, {\em On the ergodicity and mixing of max-stable processes},
  Stochastic Process. Appl., 118 (2008), pp.~1679--1705.

\bibitem{stoev2010max}
\leavevmode\vrule height 2pt depth -1.6pt width 23pt, {\em Max--stable
  processes: Representations, ergodic properties and statistical applications},
  Dependence in Probability and Statistics, Lecture Notes in Statistics 200,
  Doukhan, P., Lang, G., Surgailis, D., Teyssiere, G. (Eds.), 200 (2010),
  pp.~21--42.

\bibitem{WangOrstein}
\n{ Y.~Wang}, {\em Extremes of $q$-{O}rnstein-{U}hlenbeck processes},
  https://arxiv.org/abs/1609.00338,  (2016).

\bibitem{WangStoev}
\n{ Y.~Wang and S.~A. Stoev}, {\em On the structure and representations of
  max-stable processes}, Adv. in Appl. Probab., 42 (2010), pp.~855--877.


\end{thebibliography}
\end{document}